\documentclass{amsart}
\usepackage{amssymb,amsmath,amsfonts}

\newcommand{\C}{\mathbb C}
\newcommand{\R}{\mathbb R}
\newcommand{\N}{\mathbb N}
\newcommand{\Z}{\mathbb Z}
\newcommand{\T}{\otimes}
\newcommand{\Cn}{\C[e_0,\ldots,e_{n-1}]}
\newcommand{\an}{a_1,\ldots,a_n}
\newcommand{\Jn}{J^{\an}}
\newcommand{\Wn}{W^{\an}}
\newcommand{\ta}{\theta}
\newcommand{\Tn}{t_1,\ldots,t_n}

\newcommand{\ld}{\ldots}
\newcommand{\hk}{\hookrightarrow}

\newcommand{\Sum}{\sum\limits}
\newcommand{\Prod}{\prod\limits}

\newcommand{\nc}{\newcommand}
\nc{\g}{\mathfrak g}
\nc{\h}{\mathfrak h}
\nc{\n}{\mathfrak n_+}

\nc{\pin}{\pi_1,\ld,\pi_n}
\nc{\zn}{z_1,\ld,z_n}
\nc{\slt}{\mathfrak{sl}_2}
\nc{\Tpi}{\pi_1\T\ld\T\pi_n}
\nc{\fpi}{\pi_1*\ld*\pi_n}
\nc{\A}{\mathfrak A}
\nc{\yn}{y_1,\ld,y_n}
\nc{\pz}{\phi_Z}
\nc{\al}{\alpha}
\nc{\be}{\beta}
\nc{\CT}{(\C^2)^{\T n}}
\nc{\Ca}{\C^{a_1}\T\ld\T\C^{a_n}}
\nc{\veps}{\varepsilon}
\nc{\ik}{i_1,\ld,i_k}
\nc{\ws}{\widetilde\sigma}

\nc{\bs}{b_1,\ldots,b_s}

\nc{\qb}[2]{{\genfrac{[}{]}{0pt}{0}{#1}{#2}}_q}
\nc{\qf}[1]{(#1)_q!}

\nc{\codim}{{\mathop{\rm codim}}}
\nc{\rk}{{\mathop{\rm rk}}}
\nc{\id}{{\mathop{\rm id}}}

\newtheorem{rem}{Remark}[section]
\newtheorem{opr}{Definition}[section]
\newtheorem{theorem}{Theorem}[section]
\newtheorem{utv}{Statement}[section]
\newtheorem{lem}{Lemma}[section]
\newtheorem{cor}{Corollary}[section]
\newtheorem{prop}{Proposition}[section]

\begin{document}
\pagestyle{myheadings}
\renewcommand{\subsectionmark}[1]{}
\renewcommand{\sectionmark}[1]{}

\markboth{B.Feigin and ….Feigin}{Q-characters of the tensor products}
\title[$Q$-characters of the tensor products]
{$Q$-characters of the tensor products in $\slt$-case}
\author{B.Feigin and ….Feigin}
\address{B.F.: Landau Institute for Theoretical Physics, Chernogolovka,
142432, Russia.}
\email {feigin@mccme.ru }
\address{….F.: Independent University of Moscow, Bol'shoi Vlas'evskii per. 7,
Russia, Moscow.}
\email {evgfeig@mccme.ru }
\date{}

\begin{abstract}
Let $\pin$ be an irreducible finite-dimensional $\slt$-modules. Using the
theory of the
representations of the current algebras, we introduce a several ways to
construct a $q$-grading on $\Tpi$. We study the corresponding graded modules
and prove, that they are essentially the same.
\end{abstract}
\maketitle

\section{Introduction}

In this paper we develop a geometrical approach to the theory of
$q$-characters of the tensor products. Let $\g$ be a semisimple Lie algebra,
$\h$-Cartan subalgebra of $\g$, $\dim\h=k$. Let $\pin$ be an irreducible
$\g$-modules. There is a natural decomposition of $\Tpi$:
$$\Tpi=\bigoplus_{\gamma\in\Gamma} (\Tpi)_{\gamma},$$
where $\Gamma$ is a weight lattice. Usual character is a generating function
of the numbers $\dim (\Tpi)_{\gamma}$:
$$ch(z_1,\ld,z_k,\Tpi)=\sum_{\gamma} z_1^{\gamma_1}\ld z_k^{\gamma_k}
\dim (\Tpi)_{\gamma},$$ where
$\gamma=(\gamma_1,\ld,\gamma_k)$.
Suppose, there is some reasonable way to decompose the spaces
$(\Tpi)_{\gamma}$:
$$(\Tpi)_{\gamma}=\bigoplus_{i\in \Z} (\Tpi)_{\gamma,i}.$$
Then we can introduce the $q$-character:
$$ch(q,z_1,\ld,z_k,\Tpi)=\sum_{\gamma,i} q^iz_1^{\gamma_1}\ld z_k^{\gamma_k}
\dim (\Tpi)_{\gamma,i}.$$

There are many natural ways to introduce the $q$-characters (see \cite{sup},
\cite{mult}).
Let us briefly describe the method from \cite{fusion}.
It is based on the representation theory of the current algebras $\g\T\C[t]$.
Let $\pi$ be $\g$-module,\ $z\in \C$. Consider the homomorphism of Lie
algebras
$$\tau_z:\g\T\C[t]\to\g,\ \tau_z(g\T t^i)=z^ig,\ g\in\g, z\in\C.$$
Then we can define an action of $\g\T\C[t]$ on  $\pi$ in the following way:
$\g\T t^i\cdot v=\tau_z(\g\T t^i)\cdot v.$ We denote this so-called
evaluation representation by $\pi(z)$.

Now, let $\pin$ be an irreducible finite-dimensional representations of
$\g$, $Z=(\zn)\in \C^n$ with pairwise distinct $z_i$ (denote this
subset of $\C^n$ by $\Sigma$). In this case
$\pi_1(z_1)\T\ld\T\pi_n(z_n)$ is irreducible representation of the current
algebra $\g\T\C[t]$. Fix an arbitrary element
$u\in \pi_1(z_1)\T\ld\T\pi_n(z_n)$.
Note, that the algebra $U=U(\g\T\C[t])$ is graded by the degree of $t$:
$U=\bigoplus_{i} U_i$. Thus, we have a filtration $F_s$ on
$\pi_1(z_1)\T\ld\T\pi_n(z_n)$:
$F_s=(\bigoplus_{i=0}^s U_i)\cdot u.$ Then $\bigoplus_s F_{s+1}/F_s$
is a graded $\g\T\C[t]$-module. Denote it by $\fpi$. As $u$ is cyclic,
$$\dim(\fpi)=\dim(\pi_1(z_1)\T\ld\T\pi_n(z_n))$$
and moreover, $\fpi$ is isomorphic to $\Tpi$ as $\g$-module
($\g\subset \g\T\C[t]$). Thus, we introduced an extra grading on $\Tpi$ and
defined the corresponding $q$-character.

Our construction, a priori, depends on the choice of $Z$ and $u$. Following
\cite{fusion}, let $u$ be a product of the highest vectors of $\pi_i$ .

As $u$ is cyclic,
$$\pi_1(z_1)\T\ld\T\pi_n(z_n)\cong U(\g\T\C[t])/\bar I(Z),$$ where $\bar I(Z)$
is a left ideal in the algebra $U(\g\T\C[t])$.
Now we can formulate the following conjecture:\\
{\it Conjecture.} 1).\ The family of ideals $\bar I(Z), Z\in\Sigma$ can be
continuously extended to any point $(\zn)\in\C^n$.\\
2).\ Let $\bar I(0)$ be an extension to $(0,\ld,0)$. Then there is an
isomorphism of the $\g\T\C[t]$-modules:
$$\fpi\cong U(\g\T\C[t])/\bar I(0).$$
We will prove this conjecture in the case $\g=\slt$. (Note, that a special
case of this situation can be found in \cite{chari}).

Let $\g=\mathfrak n_-\oplus\h\oplus\n$ be Cartan decomposition with
$\mathfrak n_-$ annihilating the highest vectors of $\pi_i$.
If $Z\in\Sigma$, then $\pi_1(z_1)\T\ld\T\pi_n(z_n)$ is a cyclic
$\n\T\C[t]$-module with cyclic vector $u$.
Note, that our conjecture follows frow the similar one for the Lie algebra
$\n\T\C[t]$.

Let $\g=\slt$, $\n$ is spanned by the vector $e$. Then
$$U(\n\T\C[t])=\C[e_0,e_1,\ld],\ e_i=e\T t^i.$$
Note, that if $Z\in\Sigma$, then $\Tpi$ is cyclic as $\Cn$-module.
So, we will consider not $\bar I(Z)$, but
similary defined ideals $I(Z)\subset\Cn$. To be specific,
$$\pi_1(z_1)\T\ld\T\pi_n(z_n)\cong\Cn/I(Z)$$ as a representations of an
abelian Lie
algebra, spanned by $e_0,\ld,e_{n-1}$ (denote it by  $\A_n$). It is easy
to see, that the following elements generates $I(Z)$:
$$\bigl(\sum_{i=0}^{n-1} \al_i^k e_i\bigr)^{a_k},\ k=1,\ld,n,$$
where $a_k=\dim\pi_k$, and $\al_i^k$ are the coefficients of degree $n$
polynomial $P_k$, with a following property: $P_k(z_i)=\delta_{ik}$.

In order to study $I(Z)$, we define another family of ideals $J(Z)$,
$Z\in\Sigma$, in the ring $\Cn$ ,
with a "dual" generators:
$$\bigl(\sum_{i=0}^{n-1} z_k^i e_i\bigr)^{a_k},\ k=1,\ld,n.$$
We show, that the family $J(Z)$ can be also extended to any point
$(\zn)\in\C^n$. Connection between $J(0)$ and $I(0)$ is given by the formula:
$opp(J(0))=I(0)$, where $opp$ is $\Cn$-automorphism, $opp(e_i)=e_{n-1-i}$.

Let $\pi_1=\ld=\pi_n=\C^k$. Then it is possible to describe $J(0)$ in a
following way. Let $J^k$ be an ideal in the ring $\C[e_0,e_1,\ld]$,
generated by the coefficients of the series $e(z)^k$, where
$e(z)$ is a generating function: $e(z)=\sum_{i=0}^{\infty} e_iz^i$
(ideal $J^k$ comes from the representation theory of the Lie algebra
$\widehat\slt$). We state, that $$J(0)=J^k\cap \Cn.$$

Recall, that $I(0)$ was defined as a limit of $I(Z),\ Z\to 0,\
Z\in\Sigma$. But it is possible to define all the ideals $I(Z),\ Z\in\C^n$
at once (to be specific, not $I(Z)$, but the quotients $\Cn/I(Z)$).
Let us briefly describe the corresponding construction (call it
functional).

Let $F$ be the space of $W_1\T\ld\T W_n$-valued polynomials $f(\zn)$,
$W_i$ are linear spaces.
For all pairs $i<j$ fix a filtration on $W_i\T W_j$:
$$V^{(i,j)}(0)\subset V^{(i,j)}(1)\subset\ld\subset V^{(i,j)}
(p_{i,j})=W_i\T W_j$$
($V^{(i,j)}(k)$ are the subspaces of $W_i\T W_j$). Let
$$L^{(i,j)}=W_1\T\ld\T\widehat W_i\T\ld\T\widehat W_j\T\ld\T W_n.$$
Define a map
$\sigma^{(i,j)}:W_i\T W_j\T L^{(i,j)}\to W_1\T\ld\T W_n$:
$$\sigma_{i,j}(v_i\T v_j\T\ v_1\T\ld\T\widehat v_i\T\ld\T\widehat v_j
\T\ld\T v_n)=(v_1\T\ld\T v_i\T\ld\T v_j\T\ld\T v_n).$$
Consider $F^c$-a subspace of $F$, which consists of such polynomials $f$,
that for any pair $i<j$ and for any $k$
$$
\frac{\partial^k\! f}{\partial z_i^k}\in
\sigma^{(i,j)}(V^{(i,j)}(k)\T L^{(i,j)}), \text{ if } z_i=z_j.$$

Now, let $W_i=\C^{a_i}$, $a_i=\dim\pi_i$.
We need a special filtration on $\C^{a_i}\T\C^{a_j}$. Let
$$\C^{a_i}\T\C^{a_j}=
\bigoplus\limits_{k=0}^{\min(a_i,a_j)-1} \C^{a_i+a_j-1-2k}$$
be a decomposition of $\slt$-module on the irreducible components. Denote
$$V^{(i,j)}(l)=\bigoplus_{k=0}^{l} \C^{a_i+a_j-1-2k}.$$
Let us construct the space $F^c$ using this filtrations.
Now we define a structure of $\slt\T\C[t]$-module on $F^c$.
Note, that $F=\bigotimes_{i=1}^n(\C[z_i]\T\C^{a_i})$.
There is an action of $\slt\T\C[t]$ on $\C[z_i]\T\C^{a_i}$:
$at^k\cdot(f\T v)=(z_i^kf)\T av$. So, we have $\slt\T\C[t]$-action
on $F$. It is easy to show, that $F^c$ is submodule in $F$.
In addition, $F^c$ is a free $\C[\zn]$-module, $\C[\zn]$ acts by
multiplication.

Let $T=(\Tn)\in\C^n$, $R(T)$ is an ideal in $\C[\zn]$, generated by the
polynomials
$z_i-t_i,\ i=1,\ld,n$. Define $M(T)=F^c/R(T)F^c.$
Note, that $\slt\T\C[t]$-action on $F^c$ commutes with $\C[\zn]$-action.
So, $M(T)$ is also $\slt\T\C[t]$-module.

Recall, that $\A_n$ is a subalgebra of $\slt\T\C[t]$, spanned by
$e_0,\ld,e_{n-1}$. We will prove, that for all
$T\in\C^n$ $M(T)\cong\Cn/I(T)$ as $\A_n$-modules. Moreover,
$$M(0,\ld,0)\cong\fpi$$ as $\slt\T\C[t]$-modules.

Note, that the space $F^c$ appears in the representation theory of
$\widehat{\slt}$ as a space of the correlation functions of some vertex
operators. Special case of such situation can be found in \cite{vert}.

Our paper is organized in the following way:

In the second section we describe the ideal $J(0)$, i.e. we write down its
generators system (definition 2.2.).We prove, that
$\lim_{T\to 0} J(T)=J(0)$ (theorem 2.1.). As a consequence, we obtain a
recurrent relation for the character of $\Cn/J(0)$ (corollary 2.3.) and an
equality $J(0)=J^k\cap \Cn$ (proposition 2.3.) (note, that the recurrent
relation was obtained in \cite{fusion} by the different way).
In addition, we show, that  $I(0)=opp(J(0))$ (proposition 2.4.).

In the third section we prove the equality
$\fpi=\Cn/I(0)$ (theorem 3.1.).

The fourth section consists of the description of the functional model of
the quotient $\Cn/I(Z)$. The main statement is the theorem 4.1.

And in the last section, using the recurrent relation, we obtain the formula
(Gordon like) for
the character of the space $\Cn/J(0)$ (theorem 5.1.). One can see from the
formula, that our characters can be expressed in terms of the supernomial
coefficients from \cite{sup},\cite{mult}.

\section{Ideals and quotients}

\subsection{Definitions and main statements.}

Let $n\in \N,\ e_0,e_1,\ldots$-commuting variables.
Define generating functions $e(z)$ and $e^{(n)}(z)$ by the following way:
$$e(z)=\sum_{i=0}^{\infty}e_iz^i,\quad e^{(n)}(z)=\sum_{i=0}^{n-1}e_iz^i.$$
Let $A=(\an)\in\N^n,\ T=(\Tn)\in\C^n,$ $t_i$ are pairwise distinct.
\begin{opr}
Let $J^{(\an)}(\Tn)=J^A(T)$ be an ideal in $\Cn$, generated by the elements
$e^{(n)}(t_i)^{a_i},\ i=1,\ld,n$.
\end{opr}
\noindent Define the quotient
$$W^A(T)=W^{(\an)}(\Tn)=\Cn/J^A(T).$$
Let $\Sigma=\{T\in \C^n:\forall i\ne j\ t_i\ne t_j\}$. If $T\in\Sigma$, then
$$W^A(T)=\bigotimes_{i=1}^n \C[e^{(n)}(t_i)]/(e^{(n)}(t_i))^{a_i}.$$
Our aim is to study the limit of ideals $J^A(T), T\to 0$.

Let $p,q\in \C[\zn]$. We will write $p\div q$, if there exists
$r\in \C[\zn]$ such, that $p=qr$\ (i.e. $p$ is divisible on $q$).

Define a map $+:\Z\to \N\cup \{0\},\ s\mapsto s_+$, where
$s_+=s,$\ if $s>0, s_+=0\mbox{ otherwise}$.
\begin{opr}
$J^A$-ideal in $\Cn$, generated by the elements
$$\sum_{p_1+\cdots+p_i=s,0\le p_\al<n} e_{p_1}\ld e_{p_i}$$
for all $i\ge 0,\ 0\le s<\Sum_{p=1}^n (i+1-a_p)_+.$ In other words,
the following equality is true in the quotient $\Cn/J^A$:
$$e^{(n)}(z)^i\div z^{\Sum_{p=1}^n (i+1-a_p)_+}.$$
\end{opr}
\noindent Denote $W^A=\Cn/J^A$.
\begin{theorem}
$\lim\limits_{T\to 0} J^A(T)=J^A.$
\end{theorem}
In other words, our theorem claims, that for any way
$T(\veps)=(t_1(\veps),\ld,t_n(\veps))\in\C^n, \veps\in\R,
\lim_{\veps\to 0}t_i(\veps)=0$ we have
$\lim\limits_{\veps\to 0} J^A(T(\veps))=J^A.$

\begin{prop}
$\dim W^{(\an)}=\Prod_{l=1}^n a_l.$
\end{prop}
Note, that the proposition $2.1.$ is a corollary from the theorem $2.1.$
But this two statements will be proved simulteneously.

In our paper we also use the following $(z,q)$-bigrading of the ring
$\Cn:\ deg_ze_i=1, deg_qe_i=i$.

\subsection{Dual spaces.}
We will need the description of the dual spaces to our algebras. Let $J$ be a
homogeneous ideal (with respect to the $z$-grading) in the ring
$\Cn,\ A=\Cn/J.$
Then $A=\bigoplus\limits_k A^k$, where $k$ is $z$-degree. Let
$\ta\in (A^k)^*$.
Consider a map $(A^k)^*\to\C[z_1,\ld,z_k]$,
$$\ta\mapsto f_\ta(z_1,\ldots,z_k)=\sum_{0\le i_1,\ld,i_k<n}
z_1^{i_1}\ldots z_k^{i_k}\ta(e_{i_1}\ldots e_{i_k}).$$
Thus, $f_\ta$ is a polynomial in $k$ variables and $deg_{z_i}f<n$ for
all $i$. Note, that by this way we identify the space $A^*$ with some
subspace of
polynomials. Let us formulate two statements-descriptions of this subspace
for the ideals, defined earlier.\\
Let $T=(\Tn), A=(\an)$.
\begin{utv}
$J=J^A(T).$ Then $(W^A(T))^*$ is a sum of the spaces of polynomials
$f(z_1,\ld,z_k),\ k=0,1,\ld$, which satisfy the following conditions:\\
1). $f$ is symmetric.\\
2). $\deg_{z_j}f<n,\ j=1,\ld,k$.\\
3). $f(\underbrace{t_j,\ldots,t_j}_{a_j},z_{a_j+1},\ldots,z_k)=0$,
for all $j$:\ $a_j\le k.$
\end{utv}

\begin{utv}
$J=J^A.$ Then $(W^A)^*$ is a sum of the spaces of polynomials
$f(z_1,\ld,z_k),\ k=0,1,\ld$, which satisfy the following conditions:\\
1). $f$ is symmetric.\\
2). $\deg_{z_j}f<n,\ j=1,\ld,k$.\\
3). $f(\underbrace{z,\ldots,z}_{i},z_{i+1},\ldots,z_k)\div
z^{\sum\limits_{p=1}^n (i+1-a_p)_+},\ i=1,\ld,k.$
\end{utv}
Let us prove the following lemma:
\begin{lem}
Let $\veps\in\R, T(\veps)\in\Sigma\text{ if }\veps\ne 0 ,
\lim\limits_{\veps\to 0} T(\veps)=0$.
Let $$\ta(\veps)\in W^A(T(\veps))^*,\
\lim_{\veps\to 0}f_{\ta(\veps)}=f_\ta.$$ Then $\ta\in (W^A)^*$.
\end{lem}
\begin{proof}
Let
$$g(z,z_{i+1},\ldots,z_l)=
f_{\ta(\veps)}(\underbrace{z,\ldots,z}_{i},z_{i+1},\ldots,z_l).$$ Condition
$$f(\underbrace{t_j(\veps),\ldots,t_j(\veps)}_{a_j},z_{a_j+1},\ldots,z_l)=0$$
gives us $g\div (z-t_j(\veps))^{{(i+1-a_j)}_+}$ (all the corresponding
derivatives of $g$ are vanishing at $t_j(\veps)$). Therefore,
$g\div\prod_{j=1}^n (z-t_j(\veps))^{{(i+1-a_j)}_+}$. Now let $\veps\to 0$.
Then we will obtain the condition $3)$ from the statement $2.2$ for $f_\ta$.
\end{proof}
\begin{cor}
$\dim W^A\ge\Prod_{i=1}^n a_i$.
\end{cor}
\begin{proof}
Consider an arbitrary family $T(\veps)\in\Sigma,
\lim\limits_{\veps\to 0}t_i(\veps)=0$. Then
$$\dim\lim\limits_{\veps\to 0} (W^A(T(\veps)))^*=\Prod_{i=1}^n a_i.$$
But this limit belongs to $(W^A)^*$.
\end{proof}

\subsection{Proof of the theorem 2.1.}
We know, that $\dim (\Wn)^*\ge \prod\limits_{i=1}^n a_i$. The main goal of
this subsection is to show, that $$\dim (\Wn)^*\le \prod\limits_{i=1}^n a_i.$$
We prove this statement by the induction on $\sum a_i$.
As a consequence, we obtain the theorem $2.1.$\\
Let $a_1\le\ld\le a_{n+1}$. Denote
$$A=(\an), Aa_{n+1}=(\an,a_{n+1}).$$
\begin{lem}
Natural embedding $\Cn\hk\C[e_0,\ld,e_n]$ induces the embedding
$J^A\hk J^{Aa_{n+1}}$.
\end{lem}
\begin{proof}
It is necessary to show, that in  $\C[e_0,\ld,e_n]/J^{Aa_{n+1}}$
$$e^{(n)}(z)^i\div z^{\sum\limits_{p=1}^n (i+1-a_p)_+}.$$
In fact, $$e^{(n)}(z)^i=(e^{(n+1)}(z)-z^ne_n)^i=\sum\limits_{j=0}^i
{e^{(n+1)}(z)}^j z^{n(i-j)}\binom{i}{j} e_{n}^{i-j}(-1)^{i-j} .$$
Note, that $${e^{(n+1)}(z)}^j z^{n(i-j)}\div
z^{n(i-j)+\sum\limits_{p=1}^{n+1} (j+1-a_p)_{+}}.$$ But
$$n(i-j)+\sum\limits_{p=1}^{n+1} (j+1-a_p)_{+}\ge \sum\limits_{p=1}^{n}
(i+1-a_p)+(j+1-a_{n+1})_{+}\ge \sum\limits_{p=1}^{n} (i+1-a_p)_{+}$$ (because
$i-j+(j+1-a_p)_{+}\ge (i+1-a_p)_{+})\Rightarrow$ lemma is proved.
\end{proof}

As a consequence, we obtain a surjection
$W^A\to W_n^{Aa_{n+1}},$ where $W_n^{Aa_{n+1}}$ is a part of
$W^{Aa_{n+1}}$, generated by the action of $\Cn$ on $1$.

Suppose, that we know, that $\dim (W^A)\le\Prod_{i=1}^n a_i.$ We want to show
that $$\dim W^{Aa_{n+1}}/W^{Aa_{n+1}}_n\le
(\Prod_{i=1}^n a_i)a_{n+1}-\Prod_{i=1}^n a_i.$$ Surely, it gives us
$\dim W^{Aa_{n+1}}\le \Prod_{i=1}^{n+1} a_i$.

Denote $W=W^{Aa_{n+1}}, W_n=W^{Aa_{n+1}}_n$.
Consider a linear space $P=(W/W_n)^*$, i.¥. the space of
$f_\ta(z_1,\ld,z_l)\in W^*, \ta(x)=0\ \forall x\in W_n$.
Let us prove the following proposition:
\begin{prop}
There exists an injection $\phi:P\hk (W^{(a_1,\ld,a_n,a_{n+1}-1)})^*$.
\end{prop}
\begin{proof}
Let $f(z_1,\ld,z_{l+1})\in P$. Then $f$ doesn't contain a monomials
$$z_1^{i_1}\ld z_{l+1}^{i_{l+1}}, i_p<n \text{ for all }p=1,\ld,l+1.$$
So, $f$ can be uniquely represented as a sum
$$f(z_1,\ld,z_{l+1})=\Sum_{p=1}^{l+1}
           z_p^n g(z_1,\ld,\widehat{z_p},\ld,z_{l+1}).\leqno(*)$$
Define $\phi(f)=g$. We need to prove, that
$g(z_1,\ld,z_l)$ satisfies the condition
$$g(\underbrace{z,\ldots,z}_{k},z_{k+1},\ldots,z_l)\div
z^{\sum\limits_{p=1}^n (k+1-a_p)_{+}+(k+1-(a_{n+1}-1))_+}.\leqno (**)$$

Define
$$b_i=\#\{j:\ a_j=i\}\ (1\le j\le n+1).$$
Let $s=\max_{1\le i\le n+1} a_i$. Then $b_{>s}=0$.

Rewrite $(**)$ in terms of $b_i$:
\begin{gather}
g(\underbrace{z,\ld,z}_k,z_{k+1},\ld,z_l)\div
           z^{\Sum\limits_{p=1}^k (k+1-p)b_p},\ k<s-1,\\
g(\underbrace{z,\ld,z}_k,z_{k+1},\ld,z_l)\div
           z^{1+\Sum\limits_{p=1}^k (k+1-p)b_p},\ k\ge s-1.
\end{gather}
We will prove $(1)$ by the increasing induction from $k=1$ to $k=s-2$ and
$(2)$ by the decreasing induction from $k=l$ to $k=s-1$.

Recall, that for $f(z_1,\ld,z_{l+1})\in (W/W_n)^*$
$$f(\underbrace{z,\ld,z}_k,z_{k+1},\ld,z_{l+1})\div
           z^{\Sum\limits_{p=1}^k (k+1-p)b_p}.\leqno(3)$$
Using $(3)$ and $(*)$ one can check the correctness of $(1)$ for $k=1$, i.e.
that $$g(z_1,\ld,z_l)\div (z_1\ld z_l)^{b_1}.$$
Suppose $(1)$ is true for $k,\ k<s-1$. We have:
\begin{multline*}
f(\underbrace{z,\ld,z}_{k+1},z_{k+2},\ld,z_{l+1})=\\
=(k+1)z^n g(\underbrace{z,\ld,z}_{k},z_{k+2},\ld,z_{l+1})+
\Sum_{p=k+2}^{l+1} z_p^n g(\underbrace{z,\ld,z}_{k+1},\ld,\widehat{z_p},\ld).
\end{multline*}
We know, that
$$z^n g(\underbrace{z,\ld,z}_{k},z_{k+2},\ld,z_{l+1})\div
z^{n+\Sum\limits_{m=1}^k (k+1-m)b_m}.$$
But for $k<s-1$
$$n+\Sum\limits_{m=1}^k (k+1-m)b_m\ge \Sum\limits_{m=1}^{k+1} (k+2-m)b_m,$$
because $\sum_{i=1}^s b_i=n+1$ and $b_s>0$. Thus, for
$$H=\Sum_{p=k+2}^{l+1}
z_p^n g(\underbrace{z,\ld,z}_{k+1},\ld,\widehat{z_p},\ld)$$
we have:
$$H(z,z_{k+2},\ld,z_{l+1})\div z^{\Sum\limits_{m=1}^{k+1} (k+2-m)b_m}.$$
Let
$$h(z,z_{k+2},\ld,z_{l+1})=
        z_{l+1}^n g(\underbrace{z,\ld,z}_{k+1},z_{k+2},\ld,z_l).$$
Then
$$h=\sum_{i=0}^r z^ih_i(z_{k+2},\ld,z_{l+1}).$$
One can see, that
$$H(z,z_{k+2},\ld,z_{l+1})=\sum_{i=0}^r z^i Sym(h_i)(z_{k+2},\ld,z_{l+1}),$$
where $Sym(h_i)$ is the symmetrization of $h_i$. Thus
$$Sym(h_i)=0\ \text{ for }\ i<\Sum\limits_{m=1}^{k+1} (k+2-m)b_m.$$
One can see, that it follows, that for $h_i$ we also have
$$h_i=0\ \text{ for }\ i<\Sum\limits_{m=1}^{k+1} (k+2-m)b_m.$$
Thus
$$z_{l+1}^n g(\underbrace{z,\ld,z}_{k+1},z_{k+2},\ld,z_l)\div
z^{\Sum\limits_{m=1}^{k+1} (k+2-m)b_m}.$$
Relation $(1)$ is proved.

Now let us prove $(2)$. Suppose $l+1\ge s$. We have
$$f(\underbrace{z,\ld,z}_{l+1})=(l+1)z^n g(\underbrace{z,\ld,z}_l).$$
So $$z^n g(\underbrace{z,\ld,z}_l)\div
   z^{\Sum_{m=1}^{l+1} (l+2-m)b_m}.$$
But $l+1\ge s$, so $\sum_{i=1}^{l+1} b_i=n+1$ and
$$\Sum_{m=1}^{l+1} (l+2-m)b_m-n=1+\Sum_{m=1}^l (l+1-m)b_m.$$
Thus
$$g(\underbrace{z,\ld,z}_l)\div z^{1+\Sum_{m=1}^l (l+1-m)b_m}.$$
Now, let $k\ge s-1$ and suppose we know $(2)$ for $k+1$.
Using the relation
\begin{multline*}
f(\underbrace{z,\ld,z}_{k+1},z_{k+2},\ld,z_{l+1})=\\
=(k+1)z^n g(\underbrace{z,\ld,z}_{k},z_{k+2},\ld,z_{l+1})+
\Sum_{p=k+2}^{l+1}
z_p^n g(\underbrace{z,\ld,z}_{k+1},\ld,\widehat{z_p},\ld)
\end{multline*}
we obtain:
$$z^n g(\underbrace{z,\ld,z}_{k},z_{k+2},\ld,z_{l+1})\div
z^{\Sum\limits_{m=1}^{k+1} (k+2-m)b_m}.$$
As above, using the condition $k+1\ge s$ we obtain $(2)$ for $k$.

Thus, we proved $(1)$ and $(2)$. It gives us $(**)$ for an arbitrary $k$.
Proposition is proved.
\end{proof}
\begin{rem}
In fact, the same arguments (as in above proof) gives us, that our
injection is an isomorphism.
\end{rem}

{\bf Proof of the proposition 2.1.} $(\dim\Wn=\prod\limits_{i=1}^n a_i)$\\
From the proposition $2.2.$ $$\dim\Wn\le\dim W^{a_1,\ld,a_{n-1}}+
\dim W^{a_1,\ld,a_{n-1},a_n-1}.$$ Using the induction on $\sum a_i$
(taking into account correctness of the proposition $2.1.$ for
$a_i=1, i=1,\ld,n$), we
obtain $$\dim\Wn\le\prod\limits_{i=1}^n a_i.$$ Proposition is proved.
\begin{cor}
$\lim\limits_{T\to 0}\Jn(T)=\Jn$, i.¥. theorem $2.1.$ is correct.
\end{cor}
Since $J^A$ is homogeneous with respect to the $(q,z)$ bigrading,
$W^A$ is bigraded. Define
$$ch(\an,q,z)=\sum\limits_{k,s\ge 0} \dim (\Wn)^{k,s}z^kq^s.$$
\begin{cor}
$$ch(\an,q,z)=ch(a_1,\ld,a_{n-1},q,z)+zq^{n-1}ch(a_1,\ld,a_{n-1},a_n-1,q,z).$$
\end{cor}

\subsection{Case $a_1=\ld=a_n$}
Let $a_1=\ld=a_n=k, k\in\N$. Denote $J^k_n=J^{\an}$.
\begin{lem}
$J^k_n\hk J^k_{n+1} $.
\end{lem}
\begin{proof}
This is a direct consequence from lemma $2.2.$
\end{proof}
Let $J^k$ be an ideal in $\C[e_0,e_1,\ld]$,\ $J^k=\bigcup_{n} J^k_n$
(we regard $\Cn$ as a subalgebra of $\C[e_0,e_1,\ld]$).
\begin{prop}
$J^k$ is generated by the elements
$$\sum_{i_1+\cdots+i_k=s,i_j\ge 0} e_{i_1}\ld e_{i_k},\ s=0,1,\ld.$$
In other words $J^k$ is generated by the coefficients of the series $e(z)^k$.
\end{prop}
\begin{proof}
Let $J(k)$ be an ideal, generated by the coefficients of $e(z)^k$.
Our goal is to show, that $J^k=J(k)$. Note, that $J(k)\subset J^k$, since
$$e(z)^k-e^{(n)}(z)^k\div z^{n-1}.$$ Let $c=-\Sum_{i=0}^{\infty} z^ie_{n+i}$.
If $e(z)^k=0$, then for $i\ge k$ we have:
$$e^{(n)}(z)^i=(e(z)+z^nc)^i=\sum\limits_{j=0}^{k-1} \binom{i}{j}
e(z)^jz^{n(i-j)}c^{i-j}.$$ Thus, we obtain
$e^{(n)}(z)^k\div z^{n(i-k+1)}$ as a consequence of $e(z)^k=0$. It means, that
$J^k\subset J(k)$.
\end{proof}

\subsection{One more family of ideals.}
Consider an algebra $B$:
$$B=\bigotimes_{i=1}^n \C[y_i]/(y_i^{a_i}),\ a_i\in\N.$$
Let $Z=(\zn)\in\C^n$.
Consider $e_k\in B, e_k=\sum\limits_{i=1}^n z_i^k y_i.$
It is clear, that if $Z\in\Sigma$, then $e_0,\ld,e_{n-1}$ generates $B$. So
$$B=\Cn/I(Z),$$ where $I(Z)$ is an ideal in the ring $\Cn$. We want to
study an ideal $I(0)=I(0,\ld,0)=\lim\limits_{Z\to 0} I(Z)$.
Define the homomorphism $$opp:\Cn\to\Cn,\ e_i\mapsto e_{n-1-i}.$$
For $M\subset\Cn$ denote the image $opp(M)$ as $M^{opp}$.
Let $J(0)=J^{a_1,\ld,a_n}$.
\begin{prop}
There exists $I(0)=\lim\limits_{Z\to 0} I(Z)$ and $I(0)=J(0)^{opp}.$
\end{prop}
\begin{proof}
Recall, that $J(0)=\lim\limits_{Z\to 0} J(Z)$, where $J(Z)$
is generated by the elements
$$(\sum\limits_{i=0}^{n-1} z_k^ie_i)^{a_k}, k=1,\ld,n.$$
To prove our proposition, we construct the family of the isomorphisms
$$\phi_Z:\Cn\to \Cn,\ Z\in\Sigma$$ with a property
$\phi_Z(I(Z))=J(Z)^{opp}$
and $\lim\limits_{Z\to 0}\phi_Z=\id.$ Let
$$\pz(y_\al)=(\sum\limits_{i=0}^{n-1} z_\al^i e_{n-1-i})s_\al,\quad
s_\al=\frac{1}{\prod\limits_{\be\ne\al} (z_\be-z_\al)}$$
(if $z_i$ are pairwise distinct, then the transition matrix from $e_i$ to
$y_j$ is invertable, so it is enough to define $\pz$ on $y_{\al}$).
One can see, that $\pz$ identifies $I(Z)$ and $J(Z)^{opp}$. We have:
$$\pz (e_i)=\sum\limits_{\al=1}^{n} \pz (z_\al^i y_\al)=\sum_{j=0}^{n-1}
   e_j(\sum\limits_{\al=1}^{n} z_\al^{n-1+i-j}s_\al).$$
\begin{lem}
Let $\rho(l)=\sum\limits_{\al=1}^{n} z_\al^l s_\al$. Then $\rho(l)$
equals:\\
1). $0,\quad l=0,\ld,{n-2}.$\\
2). $1,\quad l={n-1}.$\\
3). $p(\zn),\quad l>{n-1}$, where $p(\zn)$ is a homogeneous polynomial of a
positive degree.
\end{lem}
\begin{proof}
Rational function $\rho(l)$ can be rewritten in a form:
$$\rho(l)=\frac{q(\zn)}{\prod\limits_{\al<\be}(z_\al-z_\be)},$$ where
$q$ is a polynomial. Note, that $\rho(l)$ is a symmetric function. So,
$q$ is a skewsymmetric polynomial:
$$q(\zn)=\prod\limits_{\al<\be}(z_\al-z_\be)r(\zn).$$ To obtain 1), it
is enough to mention, that $\deg r=l-(n-1)$. If $l=n-1$ then $r=const$. One
can prove, that $r=1$. 3) is obvious.
\end{proof}
Thus,
$$\pz(e_i)=e_i+\sum\limits_{j=0}^{i-1} p_j(\zn)e_j,$$ where
$p_j$ are polynomials, $\deg p_j=i-j.\Rightarrow \lim\limits_{Z\to 0}\pz=\id.$
Proposition is proved.
\end{proof}
\subsection{Extension to an arbitrary point.}
Recall, that we extended an easy defined family of ideals $I(Z), Z\in\Sigma$
to the point $(0,\ld,0)$. As a corollary, we can also describe the continuous
extension to an arbitrary point $(\zn)\in\C^n$ (in some sense, $(0,\ld,0)$
is the most "exeptional" point). The description of our extension is given
in the following statements.
\begin{utv}
Let $c\in\C,\ I(c^n)=I({\underbrace{c,\ld,c}_{n}})=
\lim\limits_{z_i\to c} I(Z)$
(as usually, $z_i$ are pairwise distinct).
Then there is an isomorphism of algebras $$\psi:\Cn/I(c^n)\to\Cn/I(0),$$
defined by $\psi(e_i)=\sum\limits_{j=0}^i c^j\binom{i}{j}e_{i-j}$ ($\psi$ is
a shift).
\end{utv}
\noindent Denote $W_n(c)=\Cn/I(c^n)$.
\begin{utv}
Let $Z=(\zn)\in\C^n$,
$$z_1=\ld=z_{\al_1}=t_1,\ld,z_{n-\al_k+1}=\ld=z_{n}=t_k,$$
$t_1,\ld,t_k$ are pairwise distinct, $\al_i$ is the number of appearance of
$t_i$ in $Z$.
Define $I(Z)$ as a limit of $I(\widetilde Z)$, all $\widetilde z_i$ are
pairwise distinct, $\widetilde z_i\to z_i$. Then $I(Z)$ is well defined and
$$\Cn/I(Z)\cong W_{\al_1}(t_1)\T\ld\T W_{\al_k}(t_k).$$
\end{utv}

Thus, we have constructed two family of ideals, depending on the point from
$\C^n$, or (in other words) two $\Prod_{i=1}^n a_i$-dimensional bundles of
factoralgebras on $\C^n$.

\section{Representations of an abelian Lie algebras and fusion product}

Recall the construction of a fusion product from \cite{fusion}.
Let $\pi$ be a representation of $\slt$. For $z\in\C$ define an action of
$\slt\T\C[t]$ on $\pi$:
$$at^k\cdot v=z^ka\cdot v\quad \forall a\in\slt, v\in\pi.$$
This $\slt\T\C[t]$-module is called an evaluation representation.
Denote it by $\pi(z)$.
Now, let $\pin$ be an irreducible representations of $\slt$,\ $\dim\pi_i=a_i$.
Let $\A$ be an abelian subalgebra of $\slt\T\C[t]$ with a base
$e_0,e_1,\ld (e_i=e\otimes t^i)$. One can check, that for $Z=(\zn)\in\Sigma$
$$\pi_1(z_1)\T\ld\T\pi_n(z_n)$$ is a cyclic $\A$-module with a cyclic vector
$u=u_1\T\ld \T u_n$ ($u_i$ is a highest vector of $\pi_i$). Introduce an
increasing filtration $F_s$ on $\bigotimes_{i} \pi_i(z_i)$:
$$F_s=<e_{i_1}\ld e_{i_k}u,\ k=0,1,\ld, \sum\limits_{j=1}^k i_j\le s>.$$
\begin{opr}
Graded $\A$-module
$$\fpi=F_0\oplus\bigoplus_{s\ge 0} F_{s+1}/F_s$$ is called a fusion product of
$\pi_1,\ld,\pi_n$.
\end{opr}
One can see, that in fact fusion product is graded $\slt\T\C[t]$-module.

Let $\A_n$ be a subalgebra of $\A$ with a base $e_0,\ld,e_{n-1}$. We will prove the
following theorem:
\begin{theorem}
$\fpi\cong\Cn/I(0)$ as $\A_n$-modules.
\end{theorem}
Note, that we used a point $Z\in\C^n$ to construct the fusion product. The
important consequence from the theorem is a fact, that the structure of
$\A_n$-module on $\fpi$ doesn't depend on $Z$. (One can obtain from the
theorem $3.1.$, that the structure of $\fpi$ as $\slt\T\C[t]$-module doesn't
depend on $Z$).

To connect $I(0)$ and the fusion product, let us redefine the latter in terms
of commutative algebra. Let $y_1,\ld,y_n$ be a commuting variables,
$y_i^{a_i}=0,\ a_i\in\N$. Consider an algebra
$$\widetilde B=\bigotimes_{i=1}^n \C[y_i]/y_i^{a_i}.$$
For $Z\in\Sigma$ define an action of $e_i, i=0,\ld,n-1$ on $\widetilde B$ by
multiplication on $\sum_{\al=1}^n z_{\al}^iy_{\al}.$
Then, $\widetilde B\cong \pi_1(z_1)\T\ld\T\pi_n(z_n)$ as a $\Cn$-module.
Thus, we obtain a filtration $F_s$ on $\widetilde B$. One can see, that it is
compatible with a structure of an algebra. Consider an adjoint graded algebra
$B$. Our aim is to prove, that $B\cong \Cn/I(0)$.
\begin{lem}
Let $W=\Cn/I$, where $I$ is homogeneous with a respect to the $z$-grading.
Introduce a filtration $F_s$ on $W$:
$$F_s=<e_{i_1}\ld e_{i_k},\ k=0,1,\ld, \sum\limits_{j=1}^k i_j\le s>.$$
Then $GrW=\Cn/I^{up}$, where
$$I^{up}=<a:\exists b\in I, b=b_0+\cdots+b_s, \deg_q b_l=l, b_s=a>.$$
\end{lem}
The proof is obvious.

We described the changing of the ideal while transfering to an adjoint
object: $I\rightsquigarrow I^{up}$. Now we will prove the following lemma
in order to connect $I(Z)^{up}$ and $\lim_{Z\to 0} I(Z)$.
\begin{lem}
Let $S_t$ be oneparameter group of homomorphisms
$$S_t:\Cn\to\Cn,\quad S_t(e_i)=t^ie_i, t\in\C^*.$$
Then:\\
$1).\quad\lim\limits_{t\to\infty}S_tI(Z)=I(Z)^{up}\\
2).\quad S_tI(Z)=I(\frac{Z}{t})=I(\frac{z_1}{t},\ld,\frac{z_n}{t})$
\end{lem}
\begin{proof}
At first, let us prove $2).$ Recall, that
$$e_i=\sum\limits_{\al=1}^n z_\al^i y_\al\Rightarrow$$
$$S_t(e_i)=\sum\limits_{\al=1}^n z_\al^i S_t(y_\al) \Leftrightarrow
t^ie_i=\sum\limits_{\al=1}^n z_\al^i S_t(y_\al)\Leftrightarrow
e_i=\sum\limits_{\al=1}^n {\Bigl(\frac{z_\al}{t}\Bigr)}^i S_t(y_\al).$$
To obtain our statement, it is enough to denote $\tilde y_\al=S_t(y_\al)$.

Now we will prove $1).$ Let $b\in I(Z), b=b_0+\cdots+b_s, \deg_q b_i=i$.
Then $$S_t\Bigl(\frac{b}{t^s}\Bigr)=\sum\limits_{k=0}^s t^{k-s}b_k.$$
We have $$\lim\limits_{t\to\infty}
S_t\Bigl(\frac{b}{t^s}\Bigr)=b_s.$$
So $I(Z)^{up}\hk\lim\limits_{t\to\infty}S_tI(Z)$ (note, that the existence
of the limit was proved earlier). But the dimensions of the quotients
are equal. Lemma is proved.
\end{proof}
{\bf Proof of the theorem 3.1.} \\
We obtain the theorem as a consequence of the lemmas 3.1. and 3.2.
$\square$

\section{Functional construction of the fusion product}
\subsection{The space of the polynomials with a conditions on the diagonals.}
Denote by $F$ the space of the $W_1\T\ld\T W_n$-valued polynomials $f(\zn)$,
$W_i$ are linear spaces.
For all pairs $i<j$ fix a filtration on the tensor product $W_i\T W_j$:
$$V^{(i,j)}(0)\subset V^{(i,j)}(1)\subset\ld\subset V^{(i,j)}
(p_{i,j})=W_i\T W_j$$
($V^{(i,j)}(k)$ are the subspaces of $W_i\T W_j$). Let
$$L^{(i,j)}=W_1\T\ld\T\widehat W_i\T\ld\T\widehat W_j\T\ld\T W_n.$$
Define a map
$\sigma^{(i,j)}:W_i\T W_j\T L^{(i,j)}\to W_1\T\ld\T W_n$:
$$\sigma_{i,j}(v_i\T v_j\T\ v_1\T\ld\T\widehat v_i\T\ld\T\widehat v_j
\T\ld\T v_n)=(v_1\T\ld\T v_i\T\ld\T v_j\T\ld\T v_n).$$
Let $f\in F, s=z_i-z_j, t=z_i+z_j.$ Consider the decomposition:
$$f(\zn)=\sum_k s^k f_k^{(i,j)}(t,z_1,\ld,\widehat z_i,\ld,
\widehat z_j,\ld, z_n).$$
\begin{opr}
$F^c$ is a subspace of $F$, consisting of such $f$, that for all $i<j$
$$
f^{(i,j)}_k (t,z_1,\ld,\widehat z_i,\ld,\widehat z_j,\ld, z_n)\in
\sigma^{(i,j)}(V^{(i,j)}(k)\T L^{(i,j)})$$ for all values of arguments.
Note, that we regard $V^{(i,j)}(k)=W_i\T W_j$ if $k>p_{i,j}$.
\end{opr}
Let all $W_i$ be finite-dimensional, $W_i=\C^{a_i}$.
Note, that $\C[\zn]$ acts on $F^c$ by multiplication.
\begin{lem}
$F^c$ is a free $\C[\zn]$-module, $\rk F^c=\prod_{i=1}^n a_i.$
\end{lem}
\begin{proof}
Firstly, $F$ is a free $\C[\zn]$-module and $\rk F=\prod_{i=1}^n a_i$
($F$ is generated by the constant polynomials, equals to the base vectors
of $\Ca$). Secondly, we can obtain $F^c$ from $F$ by the multiple iteration
of the following construction.
Let $W$ be a linear space equipped with a filtration:
$$V(0)\subset V(1)\ld\subset V(r)=W.$$
Let $G$ be a space of $W$-valued polynomials. Let
$g\in G, s=z_1-z_2, t=z_1+z_2$. Consider the decomposition:
$$g(\zn)=\sum_{i=0}^m s^ig_i(t,z_3,\ld,z_n).\eqno(1)$$
Define a subspace $G^c\subset G$:
$$G^c=\{g\in G:g_i(t,z_3,\ld,z_n)\in V(i)\T\C^{a_3}\T\ld\T\C^{a_n},\
\forall t,z_3,\ld,z_n\}.$$
Iterating this construction (substituting in general case $1,2$ by $i,j$),
we can obtain $F^c$ from $F$.
Thus, it is enough to prove, that in above construction, $G^c$ is a free
$\C[\zn]$-module.

Let $w_1,\ld,w_l$ be a base of $W$ such, that $w_1,\ld,w_{k(i)}$ form a base
of $V(i)$
($k(0)\le k(1)\le\ld\le k(r)=l$). From the decomposition $(1)$ one can see,
that $$w_1,\ld,w_{k(0)},sw_{k(0)+1},\ld,sw_{k(1)},\ld,
s^r w_{k(r-1)+1},\ld,s^r w_l$$ form a base of $G^c$ as $\C[\zn]$-module. To
finish the proof, note, that the dimension doesn't change as we replace $G$
by $G^c$.
\end{proof}
\subsection{Construction of the fusion product.}
In this subsection $W_i$ are irreducible $\slt$-modules,
$W_i=\C^{a_i}$.
Let us describe a specific filtration to be used. Let
$$\C^{a_i}\T\C^{a_j}=\bigoplus_{k=0}^{\min(a_i,a_j)-1} \C^{a_i+a_j-1-2k}$$
be a decomposition of $\slt$-module on the irreducible components. Denote
$$V^{(i,j)}(l)=\bigoplus_{k=0}^{l} \C^{a_i+a_j-1-2k}.$$
Construct $F^c$, using the above filtrations. Now, we will define an action of
$\slt\T\C[t]$ on $F^c$. Note, that $F=\bigotimes_{i=1}^n(\C[z_i]\T\C^{a_i})$.
An $\slt\T\C[t]$-action on $\C[z_i]\T\C^{a_i}$ is the following:
$at^k\cdot(f\T v)=(z_i^kf)\T av$. So, we have a structure of the
$\slt\T\C[t]$-module on $F$. One can check, that $F^c$ is a submodule of $F$.

Let $T=(\Tn)\in\C^n$, $R(T)$ be an ideal in $\C[\zn]$, generated by the
polynomials $z_i-t_i,\ i=1,\ld,n$.
\begin{opr}
$M(T)=F^c/R(T)F^c.$
\end{opr}
Note, that $\slt\T\C[t]$-action commutes with a multiplication
on the polynomials. Thus, we have a structure of $\slt\T\C[t]$-module on
$M(T)$.

Let $u_i$ be the highest vector of $\C^{a_i}$,
$\tilde u=\bigotimes_{i=1}^n u_i$. Let $u\in F^c$ be a constant polynomial,
equals $\tilde u$. Recall, that $\A_n$ is a subalgebra of $\slt\T\C[t]$
with a base $e_0,\ld,e_{n-1}$.
\begin{theorem}
For all $T\in\C^n\ M(T)\cong\Cn/I(T)$ as $\A_n$-modules. Here vector
$u\in M(T)$ is corresponding to $1\in\Cn/I(T)$.
\end{theorem}

First, we will construct an isomorphism $\Cn/I(T)\to M(T)$ for $T\in\Sigma$,
then for such $T$, that there exists exactly one pair $i\ne j$, but $t_i=t_j$.
And in the end, we will prove our theorem in its whole generality.
\begin{lem}
Theorem holds for $T\in\Sigma$.
\end{lem}
\begin{proof}
Note, that if $t_i$ are pairwise distinct, then $$\Cn/I(T)=\Ca$$ as
$\A_n$-modules (see the begining of $2.5.$). Consider the homomorphism of
$\A_n$-modules $$\psi :F^c\to \Ca,\ \psi f=f(\Tn).$$
One can see, that $\psi(R(T)F^c)=0$. So, we obtain the homomorphism
$\widetilde \psi :M(T)\to \Ca$. Let us show, that it is an isomorphism.
Really, $\Ca$ is a cyclic module and $\widetilde \psi(u)$ is a cyclic vector.
So $\widetilde \psi$ is a surjection.
In the same time $\dim M(T)=\dim(\Ca)=\prod_{i=1}^n a_i.$ Lemma is proved.
\end{proof}

\begin{utv}
For all $T\in\C^n$ there exists a homomorphism of $\A_n$-modules
$$\phi_T:\Cn/I(T)\to M(T),$$ mapping $1$ to $u$.
\end{utv}
\begin{proof}
We must check, that an action of $I(T)$ vanishes vector $u\in M(T)$.
It is obvious for $T\in\Sigma$, and because of continuity, it also holds in
the general case.
\end{proof}

\begin{lem}
Theorem $4.1.$ is true if $n=2$.
\end{lem}
\begin{proof}
Let $T=(t_1,t_2)$. We can put $t_1=t_2$. Let $t_1=t_2=0$ (general case is
an obvious consequence). Note, that it is enough to show, that $M(0,0)$ is
cyclic with a cyclic vector $u$ (see statement 4.1.). So, let
$1,y_1,\ld,y_1^{a_1-1}$ be a base of $\C^{a_1}$, $1,y_2,\ld,y_2^{a_2-1}$-a
base of $\C^{a_2}$, $e(y_j^i)=y_j^{i+1},\ j=1,2.$ In this notations $e_0,e_1$
acts on $F^c$ by the multiplication on $y_1+y_2, z_1y_1+z_2y_2$.
Note, that $2(z_1y_1+z_2y_2)=t(y_1+y_2)+s(y_1-y_2)$ (recall, that
$t=z_1+z_2, s=z_1-z_2$).
Thus, $2e_1=s(y_1-y_2)$ in $M(0,0).$

Consider a decomposition of $\C^{a_1}\T\C^{a_2}$ on the irreducible
$\slt$-components:
$$\C^{a_1}\T\C^{a_2}=V_0\oplus\ld\oplus V_r,\ \dim V_i=k_i=a_1+a_2-1-2i.$$
Let $v_1^i,\ld,v_{k_i}^i$ be a base of $V_i$. Recall, that
$$s^iv^i_j,\ i=0,\ld,r,j=1,\ld,k_i$$ form a base of $M(0,0)$.

Let us return to $e_0$ and $e_1$. We know, that
$$2^me_0^le_1^m=s^m(y_1-y_2)^m(y_1+y_2)^l$$
in $M(0,0)$. But any monomial $y_1^iy_2^j,\ i<a_1,j<a_2$
(i.¥. $\C^{a_1}\T\C^{a_2}$-base) can be rewritten as a polynomial in
variables $y_1+y_2, y_1-y_2$. Now, it is easy to show, that $M(0,0)$ is cyclic.
\end{proof}

Let $T\in\C^n$,\ $t_1=t_2$, and there is no other coinciding pairs.
Denote $M_2=F^c/R(t_1,t_2), M_{n-2}=F^c/R(t_3,\ld,t_n).$
\begin{lem}
$M(T)\cong M_2\T M_{n-2}.$
\end{lem}
\begin{proof}
Denote by $\widehat{M(T)}\subset M(T)$ the submodule, generated from $u$ by
the action of $\A_n$. We will show, that there is a surjective
homomorphism of $\A_n$-modules $\widehat{M(T)}\to M_2\T M_{n-2}.$
In this case $\dim\widehat{M(T)}\ge \Prod_{i=1}^n a_i$. But
$$\dim M(T)=\Prod_{i=1}^n a_i\Rightarrow\ \widehat{M(T)}=M(T)\cong
M_2\T M_{n-2}.$$
So, let us prove the existence of the surjective homomrphism, mapping
$$u\to \bar u=(u_1\T u_2)\T(u_3\T\ld\T u_n)$$
($u_i$ is a highest vector of $\C^{a_i}$).
$M_2\T M_{n-2}$ is cyclic, so the only thing to be proved, is a fact, that
all algebraical relations on $e_i$ which holds in $\widehat {M(T)}$, are also
fulfilled in $M_2\T M_{n-2}$. But that is practically obvious.
\end{proof}

Denote by $\Sigma_1$ the set of such $T\in\C^n$, that there exists not more
than one pair $i<j,\ t_i=t_j$. We have proved the theorem $4.1.$ for
$T\in\Sigma_1$. Before starting the proof in the general case, let us
formulate one statement.
\begin{utv}
There exists such $w_1,\ld,w_l,\ w_i\in F^c$, that for any $T$ their images
$w_i^T\in M(T)$ form a base of $M(T)$.\\
There exists such $v_1,\ld,v_l,\ v_i\in \Cn$ that for any $T$ their images
$v_i^T\in \Cn/I(T)$ form a base of $\Cn/I(T)$.
\end{utv}
\begin{proof}
Obvious.
\end{proof}
{\bf Proof of the theorem 4.1.} $(\Cn/I(T)\cong M(T))$\\
Define
$$h_i^T\in M(T),\ h_i^T=\phi_T(v_i^T) \text{ (see statement 4.1.)}.$$ Let
$A(T)$ be a transition matrix from $w_i^T$ to $h_j^T$ in $M(T)$. We know, that
for $T\in\Sigma_1\ \det A(T)\ne 0.$ We will prove, that $\det A(T)$ is a
polynomial. In this case we will obtain, that there is no
$T:\ \det A(T)=0$, because $\codim (\C^n\setminus\Sigma_1)=2.$ (In particulary,
that means, that $\det A(T)=const.$)
\begin{lem}
$\det A(T)$ is a polynomial.
\end{lem}
\begin{proof}
Let $\bar w_i^T,\bar h_j^T\in F^c$ be a representatives of $w_i^T, h_j^T$
respectively (surely, we can take $\bar w_i^T=w_i$). For $T\in\Sigma$ the
condition $(h^T)=A(w^T)$ means, that
$$\bar h_i^T(\Tn)=\sum_{j=1}^l A_{ij}(\Tn)\bar w_j^T(\Tn),\ i=1,\ld,l.$$
Note, that $\bar h_i^T,\bar w_j^T$ are polynomials. So $A_{ij}$ are rational
functions. Thus $$\det A(T)=\frac{p(\Tn)}{q(\Tn)},$$ where $p,q$ are
polynomials. But $\det A(T)$ is everywhere-defined and continuous function.
So, $\det  A(T)$ is a polynomial.
\end{proof}

Thus, we obtained, that $\det A(T)$ is a polynomial. Theorem is proved.
$\square$
\subsection{Construction of a dual module $(\fpi)^*$.}
Recall, that we fixed the $\slt$-decomposition
$$\C^{a_i}\T\C^{a_j}=\bigoplus_{l=0}^{\min(a_i,a_j)-1} \C^{a_i+a_j-1-2l}=
\bigoplus_{l=0}^s W_l.$$
Introduce the following filtration:
$$V^{(i,j)}(k)=\bigoplus_{l=0}^k W_{s-l}$$
(it is "opposite" to the filtration, used above).
Construct $\widehat F^c$, using this filtration.
Define $L(T)=F^c/R(T).$
\begin{prop}
$L(T)\cong M(T)^*$ as $\slt\T\C[t]$-modules.
\end{prop}
We will construct nondegenerated bilinear $\slt\T\C[t]$-invariant form
$(,)_T$, coupling $M(T)$ and $L(T)$.\\
{\bf Construction of the form.} Note, that there exists a unique
bilinear $\slt$-invariant form on $\C^{a_i}$,
which takes a value $1$, if we will substitute the highest and the lowest
vectors. Multiplying such forms for all $i$, we will obtain a form on $\Ca$.
Denote it by $<,>$. Thus, there is a map
$\chi:F\T F\to\C[\zn],$
$$\chi(f_1\T f_2)(\zn)=<f_1(\zn),f_2(\zn)>.$$
Let $F^c$ be constructed, using the filtration from $4.2.$
Simple and direct calculation gives us the following lemma:
\begin{lem}
Let $f\in F^c, h\in \widehat F^c$. Then
$$\chi(f\T h)\div \prod_{1\le i<j\le n} (z_i-z_j)^{\min(a_i,a_j)-1}.$$
\end{lem}
\noindent Let $\bar f\in M(T), \bar h\in L(T)$, $f,h$ are their
representatives. Define $(,)_T$ as follows:
$$(\bar f,\bar h)_T=\frac{\chi(f\T h)}
{\prod_{1\le i<j\le n} (z_i-z_j)^{\min(a_i,a_j)-1}}(\Tn).$$
One can check the correctness of the definition.
The proof of the nondegeneracy can be carried out by the same scheme, as a
proof of the theorem $4.1.$: it is obvious for $T\in\Sigma$, then it can be
checked directly for $T\in\Sigma_1$, and for an arbitrary $T$ we obtain
the theorem as a consequence of the fact, that our form is
everywhere-defined and nondegenerated outside the variety of codimension $2$.

\section{The character formula}
Recall, that in the second section we obtained a recurrent formula for the
character of $W^{\an}$ (corollary 2.3.). We write it in the following way:
$$ch(\an,q,zq)=ch(a_1,\ld,a_{n-1},q,zq)+zq^nch(a_1,\ld,a_{n-1},a_n-1,q,zq)$$
(for this relation our formulas will be a little easier). Define
$(\bs)$, where $b_i=\#\{j|\ a_j=i\},\ s=\max(a_i)$. Of course
$\sum_{l=1}^s b_l=n$. Fix a notation: $$ch(\an,q,zq)=ch(\bs).$$
\begin{lem}
The following relation is true:
\begin{gather*}
ch(\bs)=
   \sum_{j=0}^{b_s} ch(b_1,\ld,b_{s-2},b_{s-1}+j)
   z^j q^{j(n-b_s+j)}\qb{b_s}{j},\\
\qb{b_s}{j}=\frac{\qf{b_s}}{\qf{j}\qf{b_s-j}},\ \qf{k}=\prod_{i=1}^k (1-q^i).
\end{gather*}
\end{lem}
\begin{proof}
Recall, that
$$\qb{m}{k}=q^k\qb{m-1}{k}+\qb{m-1}{k-1},\
       \qb{m}{k}=\qb{m-1}{k}+q^{m-k}\qb{m-1}{k-1}.$$
Let us prove our lemma by the induction upon $b_s$. Let $n=\sum_{i=1}^s b_i$.
\begin{multline*}
ch(b_1,\ld,b_s+1)=ch(\bs)+zq^{n+1}ch(b_1,\ld,b_{s-1}+1,b_s)=\\
\shoveleft
=\sum_{j=0}^{b_s} ch(b_1,\ld,b_{s-1}+j) z^j q^{j(n-b_s+j)}\qb{b_s}{j}+\\
+zq^{n+1}(\sum_{j=0}^{b_s} ch(b_1,\ld,b_{s-1}+j+1)
z^j q^{j(n+1-b_s+j)}\qb{b_s}{j})=\\
\shoveleft=ch(b_1,\ld,b_{s-1})+
\sum_{j=1}^{b_s} ch(b_1,\ld,b_{s-1}+j)\times\\
 \times(z^j q^{j(n-b_s+j)}\qb{b_s}{j}+
          zq^{n+1}z^{j-1}q^{(j-1)(n-b_s+j)}\qb{b_s}{j-1})+\\
 +ch(b_1,\ld,b_{s-1}+b_s+1)zq^{n+1}z^{b_s}q^{b_s(n+1)}=\\
=ch(b_1,\ld,b_{s-1})+\sum_{j=0}^{b_s} ch(b_1,\ld,b_{s-1}+j)
z^j q^{j(n-b_s+j)}(\qb{b_s}{j}+q^{b_s-j+1}\qb{b_s}{j-1})+\\
+ch(b_1,\ld,b_{s-1}+b_s+1)z^{b_s+1}q^{b_s(n+1)+n+1}=\\
=\sum_{j=0}^{b_s+1} ch(b_1,\ld,b_{s-1}+j)
z^j q^{j(n-b_s+j)}\qb{b_s+1}{j}.
\end{multline*}
Lemma is proved.
\end{proof}
\begin{theorem}
\begin{multline*}
ch(\bs)=\\
=\sum_{j_{s-1}=0}^{b_s}\sum_{j_{s-2}=0}^{b_{s-1}+j_{s-1}}\ld
\sum_{j_1=0}^{b_2+j_2}
         z^{\sum\limits_{l=1}^{s-1} j_l}
q^{\sum\limits_{l=1}^{s-1} j_l(b_1+\cdots +b_l+j_l)}\times\\
\times\qb{b_s}{j_{s-1}}\qb{b_{s-1}+j_{s-1}}{j_{s-2}}\ld\qb{b_2+j_2}{j_1}.
\end{multline*}
\end{theorem}
\begin{proof}
It is a consequence from the lemma 5.1.
\end{proof}

Now, we will consider the case $a_1=\ld=a_n=k$. Our goal is to write a
formula for the character of the space $$W^k=\C[e_0,e_1,\ld]/e(z)^k$$
(notations from the section $2$). Recall, that $W^k_n$ is a subalgebra of
$W^k$, generated by the action of $\Cn$ on $1$. We will obtain the character
of $W^k$ as a limit of the characters of $W^k_n$ while $n\to\infty$.
\begin{theorem}
{\bf (Gordon's formula).}
$$ch(W^k,q,zq)=\sum_{N=(n_1,\ld,n_{k-1})}
z^{\sum\limits_{i=1}^{k-1} in_i} q^{B(N,N)}
\frac{1}{\prod_{i=1}^{k-1} \qf{n_i}},$$ where $B$ is a bilinear form with a
matrix $$B_{i,j}=\min(i,j).$$
\end{theorem}
\begin{rem}
$B(N,N)=\sum\limits_{l=1}^{k-1} (n_{k-1}+\cdots+n_{k-l})^2$.
\end{rem}
\begin{proof}
By the theorem 5.1. we have:
\begin{multline*}
ch(W^k_n,q,zq)= \sum_{j_{k-1}=0}^{n}\sum_{j_{k-2}=0}^{j_{k-1}}\ld
\sum_{j_1=0}^{j_2}z^{\sum\limits_{l=1}^{k-1} j_l}
q^{\sum\limits_{l=1}^{k-1} j_l^2}
\qb{n}{j_{k-1}}\qb{j_{k-1}}{j_{k-2}}\ld\qb{j_2}{j_1}=\\
=\sum_{j_{k-1}=0}^{n}\sum_{j_{k-2}=0}^{j_{k-1}}\ld
\sum_{j_1=0}^{j_2}z^{\sum\limits_{l=1}^{k-1} j_l}
q^{\sum\limits_{l=1}^{k-1} j_l^2}
\frac{\qf{n}}{\qf{n-j_{k-1}}\qf{j_1}\prod\limits_{l=2}^{k-1}
\qf{j_l-j_{l-1}}}.
\end{multline*}
Let us change the parameters of the summation
$$i_1=j_1,\ i_l=j_l-j_{l-1},\ l=2,\ld,k-1.$$
Denote $I=(i_1,\ld,i_{k-1})$. Note, that $j_l=\sum\limits_{p=1}^l i_p$. So
$$
ch(W^k_n,q,zq)=\sum\limits_{i_1=0}^{n}\sum_{{i_2}=0}^{n-i_1}\ld
\sum\limits_{i_{k-1}=0}^{n-\sum\limits_{p=1}^{k-2} i_p}
z^{\sum\limits_{l=1}^{k-1} (k-l)i_l}q^{\widehat B(I,I)}
\frac{\qf{n}}{\qf{n-\sum\limits_{p=1}^{k-1} i_p}}
\frac{1}{\prod\limits_{l=1}^{k-1} \qf{i_l}}\eqno (1),
$$
where $\widehat B_{i,j}=B_{k-i,k-j}$.
Note, that the $z$-degree of the summands, numbered by the indexes
$$(i_1,\ld,i_{k-1}),\ \sum\limits_{p=1}^{k-1} i_p>\frac{n}{2}$$ is greater,
than
$\frac{n}{2}$. Moreover, if $\sum\limits_{p=1}^{k-1} i_p\le\frac{n}{2}$, then
$$\deg_q\Biggl(\frac{\qf{n}}{\qf{n-\sum\nolimits_{p=1}^{k-1} i_p}}-1\Biggr)
\ge\frac{n}{2}.$$
Thus, the character of $W^k$ is a limit of
$$
\sum_{i_1+\cdots+i_{k-1}\le \frac{n}{2}}
z^{\sum\limits_{l=1}^{k-1} (k-l)i_l}q^{\widehat B(I,I)}
\frac{1}{\prod\limits_{l=1}^{k-1} \qf{i_l}}
$$
while $n\to\infty$. So, we obtain
$$
ch(W^k,q,zq)=\sum\limits_{i_1=0}^{\infty}\sum_{{i_2}=0}^{\infty}\ld
\sum\limits_{i_{k-1}=0}^{\infty}
z^{\sum\limits_{l=1}^{k-1} (k-l)i_l}q^{\widehat B(I,I)}
\frac{1}{\prod\limits_{l=1}^{k-1} \qf{i_l}}.
$$
To complete the proof, redefine $i_\al=i_{k-\al},\ \al=1,\ld,k-1$.
\end{proof}


\begin{thebibliography}{99}
\bibitem{fusion}
B.Feigin, S.Loktev, On generalized Kostka polynomials and quantum
Verlinde rule, q-alg 9812093.
\bibitem{Gordon}
B.Feigin, S.Loktev, On finitization of Gordon identities,
{\it Funct. Anal. Appl.} {\bf 35}, 53-61 (2001).
\bibitem{sup}
Anne Schilling, S.Ole Warnaar, Supernomial coefficients, polynomial
identities and q-series, Ramanujan J. {\bf 2}, 459-494 (1998); q-alg 9701007.
\bibitem{mult}
S.O.Warnaar, The Andrews-Gordon identities and q-multinomial
coefficients, Comm. Math. Phys., {\bf 184}, 203-232 (1997); q-alg 9601012
\bibitem{sto}
B.Feigin, A.Stoyanovski, Quasi-particles models for the representations
of Lie algebras and geometry of flag manifold, hep-th/9308079, RIMS 942;
Functional models for the representations of current algebras and the
semi-infinite Schubert cells, {\it Funct. Anal. Appl.} {\bf 28} (1994), 55-72.
\bibitem{vert}
B.Feigin, T.Miwa, Extended vertex operator algebras and monomial bases,
q-alg 9901067.
\bibitem{chari}
V.Chari, A. Pressley, Weyl modules for classical and quantum affine
algebras, q-alg 0004174.
\end{thebibliography}
\end{document}